# Illustrations of non-Euclidean geometry in virtual reality

„Out of nothing I have created a strange new universe." (János Bolyai)

Martin Skrodzki, RIKEN iTHEMS, Wako, Japan[1]
Email: mail@ms-math-computer.science
Twitter: @msmathcomputer2

Mathematical objects are generally abstract and not very approachable. Even experts have to work long and hard in order to obtain some understanding of higher mathematical notions. Illustrations and interactive visualizations help both students and professionals to comprehend mathematical material and to work with it. This approach lends itself particularly well to geometrical objects that are ideal for visual representation. An example for this category of mathematical objects are hyperbolic geometric spaces. When Euclid lay down the foundations of mathematics, in approximately 300 BCE, his formulation of geometry reflected the surrounding space, as humans perceive it. For about two millennia, it remained unclear whether there are alternative geometric spaces that carry their own, unique mathematical properties and that do not reflect human every-day perceptions. Finally, in the early 19$^{th}$ century, several mathematicians described such geometries, which do not follow Euclid's rules and which were at first interesting solely from a pure mathematical point of view. These descriptions were not very accessible as mathematicians approached the geometries via complicated collections of formulae. Within the following decades, visualization aided the new concepts and two-dimensional versions of these illustrations even appeared in artistic works. Furthermore, certain aspects of Einstein's theory of relativity provided applications for non-Euclidean geometric spaces. With the rise of computer graphics towards the end of the twentieth century, three-dimensional illustrations became available to explore these geometries and their non-intuitive properties. However, just as the canvas confines the two-dimensional depictions, the computer monitor confines these three-dimensional visualizations. Only virtual reality recently made it possible to present immersive experiences of non-Euclidean geometries. In virtual reality, users have completely new opportunities to encounter geometric properties and effects that are not present in their surrounding Euclidean world.

The aim of this chapter is to present the aforementioned opportunities virtual reality provides for the illustration of non-Euclidean geometries. After presenting the historical background and pre-existing mathematical, artistic, and computer graphic based two- and three-dimensional illustrations of non-Euclidean spaces, the focus lies on three-dimensional hyperbolic geometries. These are special cases of non-Euclidean spaces. Here, the chapter presents two different setups, a pure, three-dimensional hyperbolic space and a mixture of a two-dimensional hyperbolic space with one additional Euclidean dimension. Both setups exhibit different mathematical phenomena. The corresponding sections introduce these phenomena alongside with their illustrations via virtual reality. Aside from these geometry-specific considerations, it is possible to make several general statements about experiences in hyperbolic virtual reality. A separate section of this chapter is devoted to collecting such statements. Finally, the closing section consists of an outlook on what to expect from upcoming virtual reality experiences including non-Euclidean geometry. It presents a state-of-the-art overview on recent developments in the field. The following glossary collects jargon that the chapter introduces in the various sections and that is collected here as a point of reference for the reader.

---

[1] This work has been supported by the German National Scholarship Foundation as well as RIKEN.

## Glossary

| | |
|---|---|
| *Anisotropy* | Absence of isotropy, the space behaves differently, depending on the current orientation. Compare first (isotropic) and fourth (anisotropic) column in figure 4. |
| *Euclidean Geometry* | Geometry described by Euclid in 300 BCE, mathematically capturing every-day observations with a set of five postulates (see parallel postulate). |
| *Geodesic* | A straight path from start to end. In Euclidean geometry, this is a line; in spherical geometry, it is a great circle. |
| *Geodesic Divergence* | Property of hyperbolic space. Two straight path that start with the same direction, but from different points drift further and further apart. |
| *Great Circle* | Straight path on a sphere. Examples are the equator or the lines of constant longitude on a globe. See figure 1 with the segment of a great circle connecting New York and Madrid as well as two intersecting great circles. |
| *Holonomy* | Property of non-Euclidean geometries. Moving in a closed circle causes the surrounding world to rotate. See figure 1, right, where a closed path rotates the orientation by 90°. |
| *Homogeneity* | The space behaves equally everywhere, independent of the current position. |
| *Hyperbolic Geometry* | A specific non-Euclidean geometry, where not one, but infinitely many parallels exist through a point off a given line. |
| *Isotropy* | The space behaves equally, independent of the current orientation. Compare first (isotropic) and fourth (anisotropic) column in figure 4. |
| *Parallax* | Displacement in the apparent position of an object viewed along two different lines of sight. |
| *Parallel Postulate* | (also parallel axiom) One of the five postulates of Euclid, requesting the existence of exactly one parallel trough a point off a given line. |
| *Spherical Geometry* | A specific non-Euclidean geometry, where no two lines can be parallel to each other, i.e. any two non-identical lines have two intersections. |
| *Stereopsis* | Perception of depth and 3-dimensional structure obtained based on visual information obtained by two eyes. |
| *Tiling* | (also tessellation) Covering of a plane using one or more geometric shapes, called tiles, with no overlaps and no gaps. Here: also filling a space with geometric shapes without overlaps or gaps. |
| *Topologically Equival.* | Two shapes that a series of deformations transforms into each other without tearing or gluing. |

## Historical Introduction

Since their publication in approximately 300 BCE, the axioms of Euclid's *The Elements* formed the foundation of mathematical geometry and its research. It is the first application of the deductive method to mathematics and, therefore, several scholars consider it the most influential textbook ever written (Merzbach and Boyer 2011, 90). In his work, Euclid deduces the entirety of the presented geometric theorems from a set of only five postulates and several additional axioms. The first four postulates are short and concise.

"1. Let it have been postulated to draw a straight-line from any point to any point.
2. And to produce a finite straight-line continuously in a straight-line.
3. And to draw a circle with any center and radius.
4. And that all right-angles are equal to one another." (Fitzpatrick 2008, 7)

The fifth postulate – commonly referred to as the *parallel axiom*[2] – stands apart from these first four as a very convoluted formulation.

"And that if a straight-line falling across two (other) straight-lines makes internal angles on the same side (of itself whose sum is) less than two right-angles, then the two (other) straight-lines, being produced to infinity, meet on that side (of the original straight-line) that the (sum of the internal angles) is less than two right-angles (and do not meet on the other side)."[3] (Fitzpatrick 2008, 7)

Because of the obvious difference in length and content of the fifth postulate when compared to the other four, many mathematicians of the ancient world have tried to prove that this fifth statement is dependent on the other four. Thus, they would be able to derive the fifth postulate from the remaining four, rendering it unnecessary. Despite many attempts, they were unsuccessful. Over the course of the following centuries the mathematical community set aside this question, only re-discovering it in mediaeval times. Throughout the centuries to come, numerous mathematicians tried to prove dependence of the fifth postulate.[4] Again, all efforts were in vain. The situation finally changed at the turn of the 19th century.

After the multitude of fruitless attempts to prove dependence of the fifth postulate, three mathematicians – in short succession and independently of each other – proved that the fifth postulate is indeed independent of the other four. They established the exact opposite result of what the mathematical community was questioning for such a long time. The first person to recognize this fact was none other than Carl Friedrich Gauss (1777–1855), by several accounts the greatest mathematician of his time. In several letters, he hinted at his results regarding the fifth postulate, but he never formally published them, nor did he publicly announce them (Wußing 2009, 146–158). The philosophy of Immanuel Kant (1724–1804) was omnipresent at this time. Kant stated that (visual) perception is the basis of geometry: "Geometrie legt die reine Anschauung des Raumes zu Grunde" (Kant 1979, 30).[5] The creation of a notion of space that differs from the natural Euclidean space surrounding us clearly violated this metaphysical claim. Thus, Gauss feared the critique of Kantians and hid his results (Gauß 1900, 200). Nevertheless, it was Gauss, who coined the term *non-Euclidean geometry*.[6]

About twenty years after Gauss' letters, in 1832, the young Hungarian János Bolyai (1802–1860) published an appendix in a mathematical textbook written by his father. It included proof for the existence of geometries that satisfy the first four Euclidean postulates, but violate the fifth one. Thus, Bolyai's appendix proved the independence of the fifth postulate. While Gauss praised the construction in private letters to friends, he sent a harsh response to János' father, stating that he knew all the published results decades prior to their publication (Wußing 2009, 146–158). This was a strong blow to the fragile János, who never fully recovered from it. At about the same time, in Kasan, Russia, Nikolai Ivanovich Lobachevsky (1792–1856) published a treatise on non-Euclidean geometry that a local journal printed in its 1829-1830 issue (Scriba and Schreiber 2010, 421). Three more publications followed in Russian, French, and German. However, only after Gauss' death and the publication of his private letters, including his positive opinions on the matter, the ideas became widespread and acknowledged (Scriba and Schreiber 2010, 418–430).

---

[2] Even though Euclid did not present it as one of his axioms, but as one of the postulates. On the level of axioms, it is possible to deduce one formulation from another using the additional axioms provided by Euclid. However, the status of a postulate additionally implies the practicability of the described operation; see (Scriba and Schreiber 2010, 54–56).

[3] Modern textbooks teach a more accessible phrasing of this postulate: "Through a given point, only one parallel can be drawn to a given straight line." (Mazur 2006, 76) Scholars commonly refer to it as *Playfair's Axiom*, of which there are several closely related formulations, see (Ackerberg-Hastings 2017).

[4] Between 1557 and 1800, no less than one hundred works on Euclid's fifth postulate are listed, see (Engel and Stäckel 1895).

[5] More generally on Kantian philosophy and non-Euclidean geometry, see (Meinecke 1906).

[6] See (Gauss 1900, 216ff) as well as (Klein 2004, 176). For a general history of non-Euclidean geometry, refer to (Rosenfeld 1988).

Subsequent to these initial results, researchers understood that they could replace Euclid's fifth postulate by other, different demands that still yield a meaningful geometry. If through a given point off a given line, instead of one unique parallel to the line as in Euclidean geometry, there exists no parallel, this results in *spherical geometry*.[7] While the flat plane is a model for two-dimensional Euclidean geometry, *spherical geometry* resides on the sphere. A fitting analogy is the globe. A straight line along the globe describes a *great circle*[8], which goes once around the globe and comes back to its starting point. For example, the meridians are such *great circles*. However, the lines of constant latitude, parallel to the equator, are not *great circles*, but merely circular paths on the globe.[9] Therefore, when flying from New York to Madrid, a plane does not follow a line of constant latitude, but takes a segment of a *great circle* that deviates from the constant latitude and dips towards the north, see figure 1, left. As this great circle represents a straight line on the globe, it is the closest connection between the two cities. This globe analogy also illustrates how *spherical geometry* violates Euclid's fifth postulate. Namely, for any two *great circles* that are not identical, there are exactly two points of intersection, see figure 1, center. Therefore, no two *great circles* can be parallel to each other and thus no parallels exist at all in *spherical geometry*, which establishes it to be a non-Euclidean geometry.

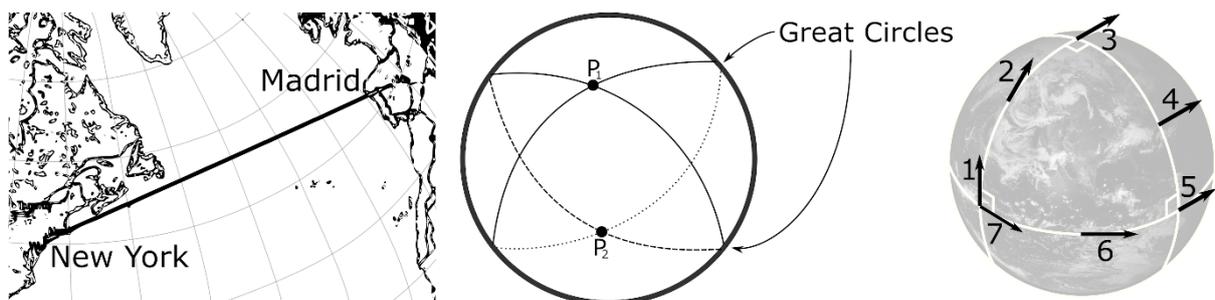

**Figure** 1 (Left: Flight Path NY – Madrid, deviating from the line of constant latitude; Center: Two non-identical *great circles* on a sphere, intersecting in two points; Right: Parallel transport of a direction along a triangle on the sphere with inner angle sum of 270°, the arrow is rotated by 90° after returning to its starting position.)

Another possibility to replace Euclid's fifth postulate leads to *hyperbolic geometry*.[10] Euclidean geometry exhibits exactly one parallel line, and there are no parallel lines in *spherical geometry*. By contrast, the demand in *hyperbolic geometry* is to have infinitely many parallels that run through a point off a given line.[11] These lines do not have to have constant distance to each other or to the original line; they merely cannot intersect the original line or each other, except for at the specified point. The two non-Euclidean geometries – *spherical* and *hyperbolic geometry* – are special cases of general classes subsumed by the broader *Riemannian geometries*.[12] They bear the name of Bernhard

---

[7] Mathematicians in Greek antiquity studied spherical geometry, but it was only much later that spherical geometry was cast within the wider framework of non-Euclidean geometries. In fact, it is necessary to change more than just the fifth postulate to obtain a meaningful axiomatization of spherical geometry, see (Whittlesey 2020), particularly chapter 3 for axiomatization.
[8] A circular path on the globe that is constructed by intersecting the globe's surface with a plane residing at the center of the globe.
[9] Except for the equator itself, which is a *great circle*.
[10] For a general discussion, see (Ramsay and Richtmyer 2013).
[11] Strictly speaking, the demand is that there are at least two parallels. However, from two parallels, by interpolating between them, infinitely many can be created. This leads to the existence of infinitely many parallels as long as the other geometrical postulates and axioms are upheld.
[12] Usually, *spherical geometry* is presented as a special case of *elliptic geometry*, which – together with *hyperbolic geometry* – can be formulated in the setting of Riemannian geometry, see (Petersen 2016).

Riemann (1826–1866), who fully realized and used the potential of non-Euclidean geometry in his works. Subsequently, Hermann Minkowski (1864–1909) developed a geometry encompassing the usual three dimensions of space, adding time as a fourth dimension (Oloff 2018, 148). Finally, it was in the language of these non-Euclidean geometries, like the *Minkowski spaces*, that Albert Einstein (1879–1955) was able to frame his general theory of relativity. One example for this framing is the *rubber sheet analogy*, an image presumably developed by Einstein himself[13] to explain his theories.[14] The rubber sheet is flat while lying on a table. However, when one lifts it up and ballasts it with an object, it curves downwards due to gravity. This two-dimensional analogy illustrates how three-dimensional space bends around heavy objects such as stars or black holes. This bent space then follows a non-Euclidean geometry, for instance light travels in curves around black holes instead of in straight lines. Aside from astrophysics, non-Euclidean geometry found numerous applications in architecture (Gawell 2013), electrical engineering (Bolinder 1958), and in a multitude of other scientific fields (Jennings 2012).

## Two-dimensional Illustrations

The relevance of non-Euclidean geometries for mathematics, physics, and other applications motivates a need for thorough introductions to this topic. Visualizations are particularly important to help students grasp the non-intuitive aspects of non-Euclidean geometries. Different researchers have developed several visualizations of hyperbolic spaces. All these bear the names of famous and accomplished mathematicians. However, they illustrate the two-dimensional hyperbolic plane as a part of the common Euclidean plane; therefore, they cannot preserve all features of hyperbolic space. The *Beltrami-Klein model*, for instance, renders the hyperbolic plane into a disk in Euclidean space, see figure 2, left. It does map straight hyperbolic lines to straight lines in the model, but it distorts the angles between these lines.[15] The namesakes are Eugenio Beltrami (1835–1900) and Felix Klein (1849–1925).

Another model, developed in the same timeframe, does preserve angles between lines, but it depicts straight lines as curved, see figure 2, center. Henri Poincaré (1854–1912) who also used a disk as embedding space developed this latter model. Literature now refers to this mapping as the *Poincaré disk model*. It is important to distinguish Poincare's disk model from a third popular model, which Poincaré also developed. The *Poincaré half-plane model* utilizes the entire upper half of the Euclidean plane to map the hyperbolic plane. It renders straight lines either as vertical lines or half circles while also preserving angles, see figure 2, right.

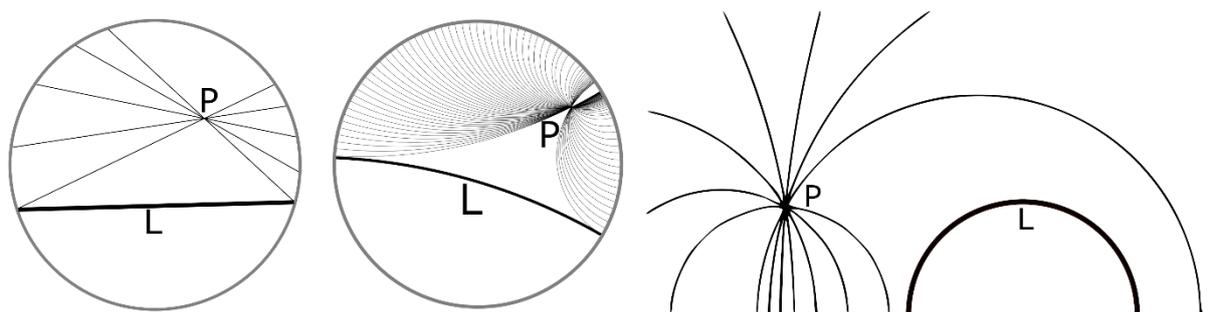

**Figure** 2 (A line $L$ and several parallels to $L$ through a point $p$ in different models of the hyperbolic plane. Left: Beltrami-Klein model in a disk, straight lines are straight; Center: Poincaré disk model, straight lines are curved; Right: Poincaré half-plane model, straight lines are vertical lines or half circles, but angles are kept intact.)

---

[13] Confer Einstein's letter from March 1917 to Willem de Sitter, as printed in (Hentschel 1998) comparing curved space to a cloth.
[14] For a general discussion of the rubber sheet analogy, see (Kersting and Steier 2018).
[15] Mathematically speaking, the map is not *conformal*.

From these listed models, the *Poincaré disk model* in particular has found widespread use. The work of Maurits Cornelis Escher (1898–1972) popularized it beyond mathematics. In his *Circle Limit* series (1958–1960), Escher drew his inspiration from a diagram, which the British mathematician Donald Coxeter (1907–2003) used in a 1957 mathematical paper on hyperbolic geometry.[16] In a series of four woodcuts, Escher combines properties of hyperbolic geometry, such as triangle angle sums larger than 180°, with his iconic interlocking tilings as well as rotational symmetries.[17] Since the release of Escher's woodcuts, hyperbolic geometry is a well-established and widely used technique in arts.[18]

Recently, software developers and artists have demonstrated the potential of combining non-Euclidean geometries and new media, like computer games and simulations. An example of this fusion is the game *Hyperrogue*[19], a two-dimensional *rogue-like*[20] turn based game, which is set in the hyperbolic plane. The two game developers, Eryk Kopczyński and Dorota Kopczyńska, have a background in mathematics and computer science research. Their game uses the properties of hyperbolic space as game mechanics that the players must use to their advantage. Thereby, it teaches an understanding of this unintuitive geometry in a playful way.

Another example is the *Curved Spaces* simulation by mathematician Jeffrey R. Weeks, which makes use of state-of-the-art computer graphics algorithms to render several non-Euclidean geometries (Weeks 2002). The applications relates the three-dimensional hyperbolic space with the three-dimensional surface of a four-dimensional ball. Within this space, Week's simulation allows the user to explore "interesting visual effects and rich symmetries not available in flat space models" (Weeks 2002, 90).

## Virtual reality illustrations of hyperbolic geometry

The natural feeling of the surrounding Euclidean geometric world hinders the experience and exploration of other geometric settings. While two-dimensional illustrations – like the *Hyperrogue* game or Escher's *Circle Limit* woodcuts – are helpful educative tools, they remain constructed artefacts on a canvas or a computer screen. As soon as the users look away from them, they are back to their regular Euclidean space. Furthermore, the two-dimensional illustrations provide an outside view of the hyperbolic space imperfectly mapped into a part of the Euclidean plane. A natural next step is to consider three-dimensional hyperbolic space. However, just as the *Curved Spaces* visualization places three-dimensional hyperbolic space onto a monitor in a surrounding Euclidean world, a *Poincaré ball model*[21] of three-dimensional hyperbolic space only provides an extrinsic view, embedded into a surrounding three-dimensional Euclidean world.

Virtual reality can eliminate this coexistence between different geometrical setups by providing a fully immersive experience of the non-Euclidean space. With more and more affordable hardware for virtual reality, further applications become available to teach hyperbolic geometry in an immersive way. In addition, virtual reality hardware allows users to walk around freely – to a certain extent – within this non-Euclidean world. It tracks the head movements of a user and allows them to kneel or lie down, providing an interactive way of experiencing these mathematical concepts. For such simulations, "[the] goal is to make three-dimensional non-euclidean [*sic*] spaces feel more natural by giving people experiences inside those spaces, including the ability to move through those spaces with their bodies"(Hart et al. 2017a, 34).

---

[16] For the paper of Coxeter, see (Coxeter 1957). For a general discussion of the *Circle Limit* series with regard to hyperbolic geometry, see (Herfort 1999).
[17] Specific mathematical aspects of the *Circle Limit* series are discussed in (Dunham, Lindgren, and Witte 1981) and (Dunham 1999).
[18] The *Journal of Mathematics and the Arts* and the archive of the mathematical art conference *Bridges* list a multitude of entries for the keyword *hyperbolic*. For a broader discussion, see (Henderson 1983).
[19] For a description of the game, consider the corresponding publication (Kopczyński, Celińska, and Ctrnáct 2017).
[20] *Rogue*-like games are a sub-category within the broader category of *role-playing games (RPG)*. Procedurally generated levels and the absence of save games that force the player to restart if the character dies characterize them.
[21] The three-dimensional generalization of the *Poincaré disk model*.

The following sections are devoted to the presentation of several corresponding projects, which visualize mathematical concepts that users cannot explore while moving in their three-dimensional Euclidean world. These projects are interesting from several perspectives. First, they provide a means for teaching non-Euclidean geometries, for instance in the university classroom. Second, they form the basis for ongoing research projects on aspects of hyperbolic space, such as visualization of large-scale data sets. Third and finally, they have artistic value and are accessible because of their immersive controls. Thus, these projects are suitable for mathematical science communication activities.

## Explorations of $\mathbb{H}^3$

A first – and arguably the easiest – setup to discuss is three-dimensional hyperbolic space, denoted by $\mathbb{H}^3$.[22] In this space, the geometry is the same at every point, just as Euclidean space behaves the same independent of the observer's position within it. This indifference to position is the property of *homogeneity*.[23] Furthermore, when picking an arbitrary point in space, the geometry behaves equally, independent of the direction one moves from this point. This mimics the user's experience of three-dimensional Euclidean space, which also behaves equally in all directions. This is the property of *isotropy*.[24] Regarding these two properties, three-dimensional hyperbolic space is equal to three-dimensional Euclidean space. It does not affect a person's step size whether they start in one corner of a room or in another corner. Similarly, one stride covers the same length independent of whether they move towards the center of the room or along the walls[25].

Given an implementation of three-dimensional hyperbolic space $\mathbb{H}^3$ in virtual reality, one of the features of the surrounding space that a user can experience is the concept of *holonomy*.[26] It describes a peculiar phenomenon: When the user moves their head in a circular motion the world slowly turns around them, even though they are not turning their head, see figure 4, first column. Any curved manifold – such as the globe – exhibits holonomy. Imagine an ant starting at the equator of the globe, right on the zero meridian, holding an arrow pointing along the meridian towards the North Pole, see figure 1, right, position 1. Holding the arrow straight, it walks along the meridian until it reaches the pole; see figure 1, right, position 3. Without turning the orientation of the arrow that it is holding, it takes steps sideways to its right until it reaches the equator; see figure 1, right, position 5. Because it always stepped sideways and never changed the orientation of the arrow, the arrow now points straight along the equator. From here, it steps backwards along the equator to the zero meridian, where it started, not changing the orientation of the arrow, which now points along the equator, see figure 1, right, position 7. Having walked this awkward triangle on the globe the ant is back where it started its walk. However, initially, its arrow was pointing to the North Pole while now it is pointing along the equator. This is the effect of holonomy: The ant moved in a circle and its orientation changed. The same general behavior holds in $\mathbb{H}^3$. Because of holonomy, there are two possible ways for a user to focus on an object behind them. Either, they turn around to look at the object, or they move their head in concentric circles, which causes the world around them to turn until the object comes into sight.

---

[22] See (Hart et al. 2017a) for a complete discussion of this virtual reality experience.
[23] From the Greek ὁμός (homos) "same" and γένος (genos) "kind".
[24] From the Greek ἴσος (isos) "equal" and τρόπος (tropos) "way".
[25] The universe is not a three-dimensional Euclidean space. The special and general theory of relativity explain that space can bend and thus exhibit phenomena where a curved path, for example around a black hole, is shorter than a straight path, for instance from the sun to earth. However, these effects are negligible in our daily observations, which caused us to develop most pre-Einsteinian physics within three-dimensional Euclidean space. See (Oloff 2018) for a discussion of space-time geometry.
[26] From the Greek ὅλος (holos) "entire" and νόμος (nomos) "law", where the latter part is also derived from the Indo-European root for "number". While the exact origin of the word is unclear, it could be related to the mathematical expression *holomorphic*, "entire shape". It is a measure of the extent to which parallel transport around closed loops fails to preserve the transported geometrical data.

The floor in a three-dimensional hyperbolic virtual reality exhibits another phenomenon of three-dimensional hyperbolic space $\mathbb{H}^3$. It is clear from the introduction that for a straight line in hyperbolic geometry and some point not on this line, there are infinitely many lines running through that point, which do not intersect the given line. While the lines all intersect in the given point, they diverge with increasing distance to the point, that is to say the distance between the lines grows while moving away from the point, see figure 2. This is the phenomenon of *geodesic*[27] *divergence*. Users in virtual reality observe a specific case of this when they look straight ahead. Their line of sight thus forms a straight line in three-dimensional hyperbolic space. Furthermore, their feet mark a point offset from this line of sight and the floor provides a parallel to the line of sight. By *geodesic divergence*, the distance of the floor to the line of sight has to grow with increasing distance from the user and therefore, the users will see the floor gradually descending downwards away from them, see figure 3, left. If they move along their line of sight, the floor will fall out from under their feet and they will start to float, see figure 3, left. While the experience of these two phenomena – *holonomy* and *geodesic divergence* – possibly causes virtual reality sickness, it also enables the user to experience these two otherwise abstract mathematical notions firsthand.

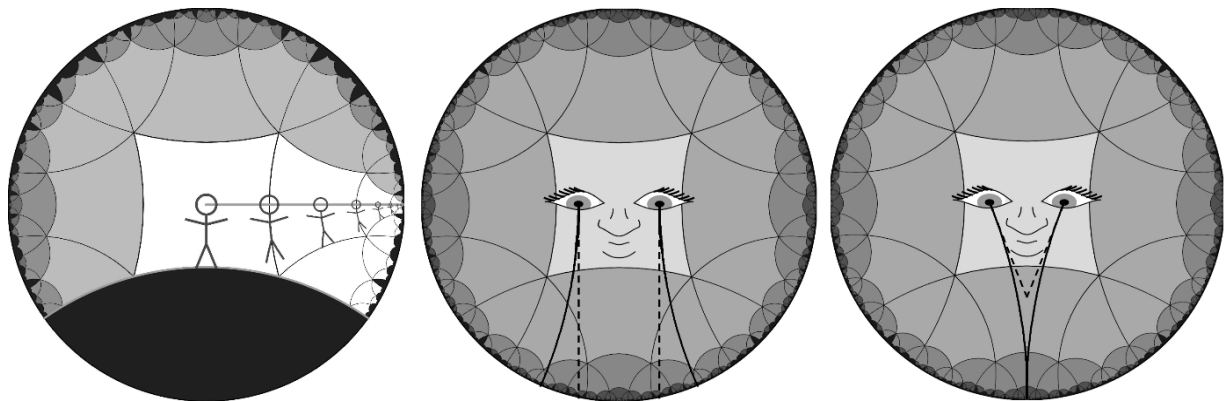

**Figure** 3 (Left: *Geodesic divergence*, the floor is a parallel line to the line of sight and thus falls away from the latter; Center: *Geodesic divergence* for eye rays; Right: Users have to squint to look at far away objects; Illustration reproduced with kind permission from (Hart et al. 2017a; Hart art al. 2017b).)

## Explorations of $\mathbb{H}^2 \times \mathbb{E}$

Both *holonomy* and *geodesic divergence* hinder the easy exploration of $\mathbb{H}^3$ in virtual reality, as it is extremely involved to move to a given point while keeping the orientation frame leveled horizontally and having the floor fall away. A combination of hyperbolic and Euclidean space alleviates these problems, as it "enables the user to traverse the hyperbolic plane horizontally as they walk through a room, yet it retains familiar euclidean [*sic*] geometry in the vertical direction" (Hart et al. 2017b, 42). This combined space consists of a two-dimensional hyperbolic space $\mathbb{H}^2$ as a plane to walk on with a one-dimensional Euclidean space $\mathbb{E}$ extending orthogonally from this plane. Similar to $\mathbb{H}^3$, this space is *homogeneous*. However, because of the mixture of the two different geometries, the space $\mathbb{H}^2 \times \mathbb{E}$ is *anisotropic*, meaning the geometry is different when looking along one of the $\mathbb{H}^2$ directions compared to the $\mathbb{E}$ direction. Compare the undistorted view of $\mathbb{H}^3$ in figure 4 (a) with the seemingly distorted views of $\mathbb{H}^2 \times \mathbb{E}$ in figures 4 (b) and (d).

---

[27] From Greek γεωδαισία (geodaisia) "division of Earth". See (Kühnel 2015, 276ff).

When navigating $\mathbb{H}^2 \times \mathbb{E}$ spaces, the missing *isotropy* causes a notable visual feature. Namely, the aspect ratio of objects change as one moves towards or away from them. Since the height of objects behave according to the Euclidean dimension, their height will grow linearly in accordance with the distance moved towards them. However, the width of objects obey the hyperbolic dimensions of the space; therefore, objects' width grow exponentially as the user moves towards them. Note that this visual experience appears in neither three-dimensional Euclidean nor three-dimensional hyperbolic space as both width and height scale equally in these two environments. Thus, a user can immediately experience *anisotropy* in a $\mathbb{H}^2 \times \mathbb{E}$ virtual environment. Compare to the undistorted octagon in figure 4 (a) and the distorted octagon in figure 4 (d).

Another noteworthy experience the user makes in this geometrical space is direction-dependent *holonomy*.[28] As discussed above, in $\mathbb{H}^3$, moving in any closed loop in space causes a rotation of the space. Because of the *anisotropy* of $\mathbb{H}^2 \times \mathbb{E}$, moving along loops does not necessarily have this effect. For instance, when following the set of movements "right, up, left, down," which cause rotation in three-dimensional hyperbolic space, see figure 4, first column, nothing happens in $\mathbb{H}^2 \times \mathbb{E}$, see figure 4, fourth column, just as expected from three-dimensional Euclidean space. This happens because movements up and down involve the Euclidean dimension of the space and thus no *holonomy* occurs.[29] However, when looking up – along the Euclidean dimension – and walking in a loop "forward, left, backward, right," this loop lies entirely in the hyperbolic dimensions and thus causes a rotation of the orientation frame comparable to that in three-dimensional hyperbolic space, see figure 4, second column.

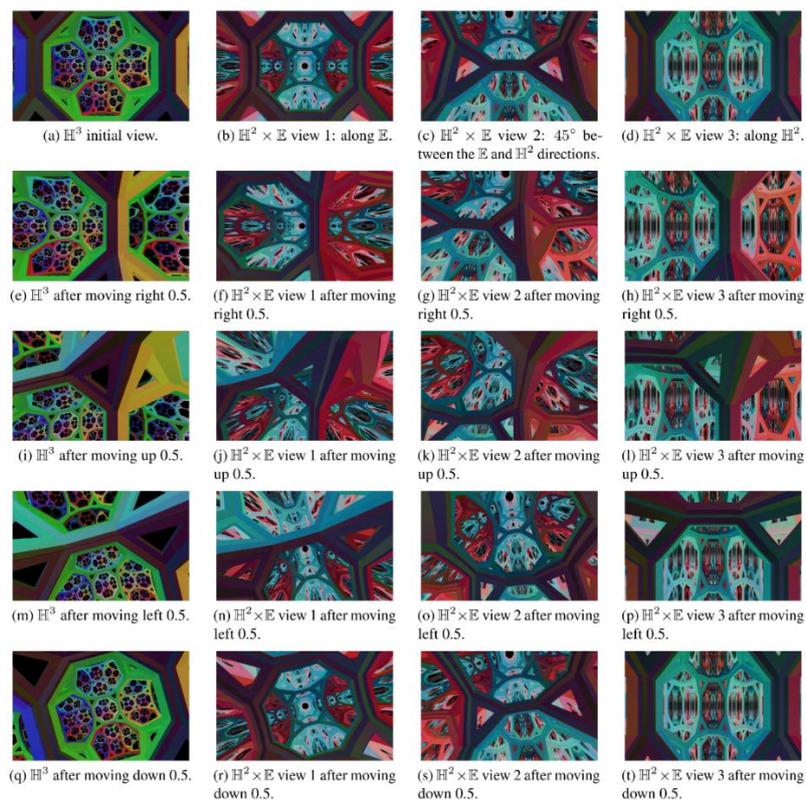

**Figure** 4, (Illustrating three-dimensional hyperbolic space as well as $\mathbb{H}^2 \times \mathbb{E}$ and movements in both spaces; Illustrations reproduced with kind permission from (Hart et al. 2017b).)

---

[28] See footnote 26.
[29] More specifically, because isometries in $\mathbb{H}^2$ and $\mathbb{E}$ commute in the Cartesian product $\mathbb{H}^2 \times \mathbb{E}$, see (Hart et al. 2017b, 45).

Because of the planar structure of the hyperbolic portion of $\mathbb{H}^2 \times \mathbb{E}$, the floor does not fall away and virtual reality sickness is not prevalent. Still, users can experience the phenomenon of *holonomy* firsthand, with a set of simple movements, which makes the non-intuitive rotation of space accessible to users without a thorough understanding of mathematical formulae.

## Experiences in hyperbolic virtual reality

The previous section discussed experiences specific to three-dimensional hyperbolic space and $\mathbb{H}^2 \times \mathbb{E}$. However, there are several effects and phenomena that users can experience in either of these geometrical realms. A first one to discuss is at the heart of modern display tools for three-dimensional perception.

Virtual reality devices, just like shutter glasses used for 3D films and even old-fashioned red/green glasses, create an illusion of a three-dimensional scene by providing the left eye of the user with a slightly different view than their right eye. The brain generates an illusory impression of depth from these differences. The general concept of such inference of three-dimensional information from binocular vision is *stereopsis*.[30] Humans employ several strategies within this realm, one of them being *parallax*[31], which uses the difference in the appearance of objects when viewed from several viewpoints. While astrophysicists employ parallax at extremely large scales, to determine the distance to other solar systems for instance, we use parallax to estimate distances to objects in everyday vision within Euclidean space. Focusing on an object at a great distance requires looking exactly straight, while objects that are just a few centimeters away require squinting. Consequently, humans can give estimates for the distance of objects by the position of their eyes, see figure 5, left.

In both hyperbolic settings presented – $\mathbb{H}^3$ and $\mathbb{H}^2 \times \mathbb{E}$ – parallel rendering rays 'shot' from two points will eventually diverge, a phenomenon called *geodesic divergence*. This phenomenon is responsible for the perception of the floor falling away from the user, as described above. When a user looks straight ahead, the information reaching either eye lies on one of two parallel lines, one for each eye respectively. However, these lines diverge in hyperbolic space, see figure 3, center. To actually focus on an object far away, the users have to point their eyes inwards, see figure 3, right. To someone who is used to viewing in Euclidean space, as discussed above, this squinting indicates closeness. The consequence of these effects is "that all objects appear to be relatively close by" (Hart et al. 2017b). That is, by placing the virtual 'eyes' of the users using the average interpupillary distance[32] as a measure, the simulation would "display the world as a hyperbolic being would see it" (Hart et al. 2017b). Experiments with reversing glasses show that humans adapt to different visual stimuli after a short period of days (Erismann and Kohler 1950). It remains a topic of open research to investigate whether users can learn to navigate hyperbolic space with this profoundly altered depth perception, which is fundamentally different to the Euclidean spaces they are familiar with.

While movement and orientation with diverging visual rays is an educational experience within a hyperbolic geometry simulation, it hinders quick and unprepared exploration of the space, as users cannot rely on their modes of depth perception learned in the Euclidean environment. As an alternative, the software renders the visual cues to the user not within the hyperbolic geometry, but within its model as part of Euclidean space. That is, the computation runs within the three-dimensional analogues of the models in figure 2. This carries two immense benefits: first, the user can utilize their common depth perception, and second, the virtual reality controllers appear at their exact position corresponding to the real, Euclidean world. However, it robs the user of the experience of perceiving the world as a hyperbolic being. Thus, both approaches have their respective educational advantages.[33]

---

[30] Derived from the Greek στερεο (stereo) "solid" and ὄψις (opsis) "appearance, sight". See (Howard and Rogers 1995) for a general discussion of the concept.
[31] Derived from the Greek παράλλαξις (parallaxis) "alternation".
[32] This interpupillary distance is $59.99 \pm 3.03$ for adult females and $62.21 \pm 3.27$ for adult males on average (MacLachlan and Howland 2002).

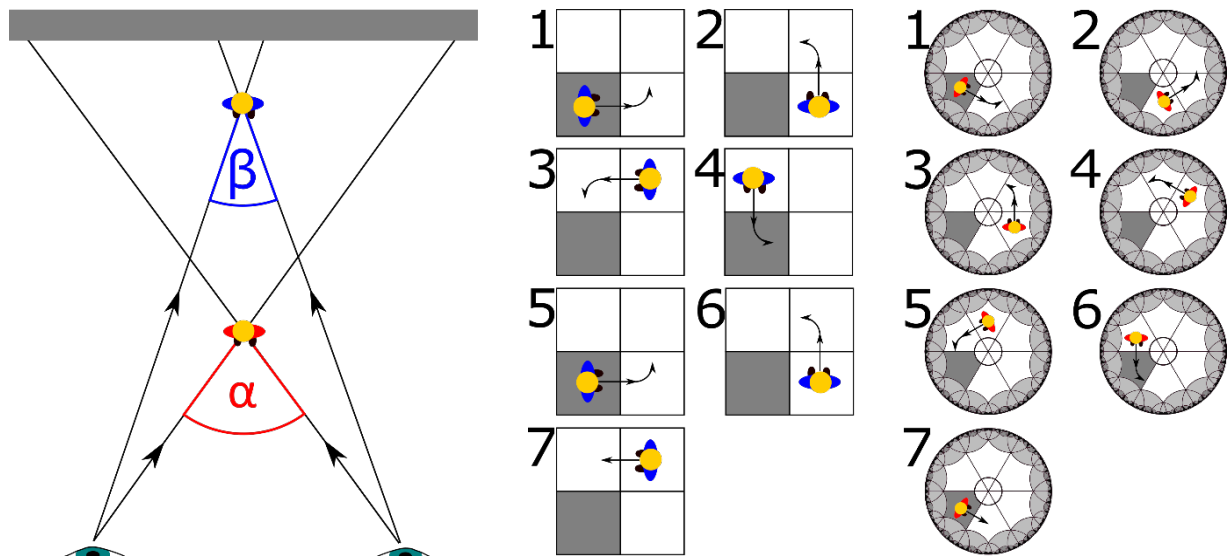

**Figure** 5, (Left: Parallax for distance estimation, the person in front, with the larger angle $\alpha$ is perceiver closer from the background than the person with the smaller angle $\beta$; Center: Parallel transport in Euclidean geometry. Here, four rotations of 90° correspond to once cycle; Right: Parallel transport in hyperbolic geometry. Here, six rotations of 90° correspond to one cycle; All Illustrations from Hart+17.)

The discussion up to this point omitted a crucial factor. Namely, in an empty hyperbolic space, users would not notice any of the described effects because they would not have any objects in the space to orient themselves by. Thus, it is vital to fill the space, ideally in a way that indicates the geometry's features by the objects that populate it. A natural way of filling Euclidean space is to stack cubes next to and on top of each other, a three-dimensional analogue to a two-dimensional checkerboard design. Thus, the user moves within a cubical room, with a square floor, a square ceiling, and four square walls. Passing through one of the walls (or floor or ceiling for that matter) leads into a neighboring cube that has the exact same dimensions. Walking around an edge formed by two walls takes the user into three different cubes before returning them to their starting point, the original cube.[34]

However, the described cubical arrangement does not tile hyperbolic space.[35] Depending on the model used, cubes can give a possible separation of the space, where not four, but six cubes fit around an edge.[36] This structure imposed on the hyperbolic spaces illustrates another hyperbolic phenomenon different from *parallax*. It relates to the inner angles of hyperbolic polygons. Due to the curved nature of hyperbolic space, the angle sum in a triangle is strictly less than 180° (Ramsay and Richtmeyr 2013, 48–49). However, to the user of the virtual reality experience, as they move through the cubes each of the floor squares will appear as a regular, Euclidean square with an angle sum of 360°.

---

[33] Technically, the effect on the perception is determined by the magnitude of curvature imposed on the virtual hyperbolic space. Compare to (Hart et al. 2017b, 47).
[34] In *Schläfli* symbol notation, this cubical tiling is written as {4,3,4}, meaning that the tiling cells are made from square faces (four sides), three of the faces meet at a vertex, and four cells are arranged around an edge.
[35] This is due to the negative curvature of hyperbolic space.
[36] That is to say, this works for specific magnitudes of negative curvature. Again, in *Schläfli* symbol notation, this cubical tiling is written as {4,3,6} as its cells are composed from squares, three of the faces meet in a vertex, but now six cells are arranged around an edge.

When observing users navigating a virtual world, an interesting method is to mark a square grid in the real world and have the user start from the same point in a real and a virtual square, see figure 5, right, position 1. If one movement consists of turning left 90° and walking straight into the center of the next square, then in hyperbolic geometry with the tiling discussed above, the users have to do six such movements in order to be back in the cube where they started, see figure 5, right, positions 2–7. From the point of view of the user within the virtual space, they have moved in a closed circle and are back at their original position. However, their physical body in the real world will not be back where they started, but will have moved two additional squares. This is a powerful physical manifestation of the underlying difference between Euclidean and hyperbolic space.[37]

## How to proceed with VR and mathematics?

The discussion so far illustrates how virtual reality environments can visualize hyperbolic geometric concepts for a user. Several phenomena of non-Euclidean spaces translate directly to virtual reality experiences. These are available for education, science communication, and research. However, this is not the end of the story, as there are further directions to explore, some of which the remainder of the chapter will briefly summarize.

### More immersive experiences by mapping hands and placing objects

As setups with virtual reality allow the user to move around freely, they also allow for an interactive exploration of *parallax* or *angle sums*. Other simulations or games provide an enhanced sense of immersion in virtual reality by mapping the user's hands or hand-held controllers into the virtual realm. Phenomena like *geodesic divergence* make this an intrinsically hard problem, as discussed above. Thus, it is impossible to combine realistic mapping of controllers together with the illusory depth perception of a hyperbolic being.

Still, being able to interact with objects in the virtual realm further increases the sense of immersion. Naturally, this would require rendering said objects in hyperbolic space. However, the majority of algorithms that render Euclidean space make use of straight lines as geodesics.[38] In hyperbolic space, these efficient renderings are not available, making the process more costly and therefore unusable for the high-framerate requirements of virtual reality. The rendering of the population of the space using cubical tiling as described above is only possible by a trick. Namely, all movement of the user takes place in one pre-rendered cube. When the user moves through a wall into another cube, they simply enter the central cube again, from the other side, like walking through a portal and returning to the same room from the opposite side. Since no objects exist that the user can pick up or leave behind, it is easy to maintain this central viewpoint. Adding such objects would call for tremendously different rendering algorithms.

### Visualization of the Thurston geometries

The hyperbolic spaces $\mathbb{H}^3$ and $\mathbb{H}^2 \times \mathbb{E}$ investigated above play a major role in the geometrization conjecture of the mathematician William Thurston (1946–2012).[39] The conjecture is set in the mathematical field of topology. This field is concerned with the classification of objects by their most fundamental shape; therefore, two objects are *topologically equivalent* if a series of operations can deform them into each other without tearing or gluing. For example, a cube and a ball are *topologically equivalent*. However, they differ from a donut, as an operation would need to glue the hole of the donut in order to deform it into a ball or a cube. Within this setting, the entire set of finite two-dimensional surfaces without boundaries is fully classified. It consists of the surface of a ball, the surface of a donut, the surface of a two-holed torus, the surface of a three-holed torus, and so on. These nicely correspond to the different two-dimensional non-Euclidean geometries. The sphere has spherical geometry, the torus has Euclidean geometry, and all multi-hole tori have hyperbolic geometry.[40]

---

[37] See (Segerman et al. 2017a) and (Segerman et al. 2017b) for further explanations and a demonstrations of the different effects discussed here.

[38] See for instance (Shirley and Morley 2003).

[39] Thurston stated the conjecture in 1980 and Grigori Jakowlewitsch Perelman (1966*) proved it in 2002. It is strongly tied to the much older *Poincaré conjecture*. See (Szpiro 2008) for a popular science introduction.

The geometrization conjecture of Thurston aims at a similar classification, but with reference to three-dimensional analogues of the objects listed above. While there is no complete characterization of the entire set, as in the two-dimensional case, it is possible to classify at least certain pieces of the three-dimensional objects. After decomposing an object into said pieces, each of these exhibits one of the following geometries: three-dimensional spherical geometry $\mathbb{S}^3$, three-dimensional Euclidean geometry $\mathbb{E}^3$, three-dimensional hyperbolic geometry $\mathbb{H}^3$, one of two mixtures of these ($\mathbb{S}^2 \times \mathbb{E}$, $\mathbb{H}^2 \times \mathbb{E}$), and three more complicated mixtures going by the names of $Nil$, $PSL_2(\mathbb{R})$, and $Solv$ (Thurston 1997). These last three also mix Euclidean spaces and hyperbolic spaces, but they include "twists" and "stretches" in the mix.

The discussion above already established that virtual reality is tremendously helpful in visualizing $\mathbb{H}^3$ and $\mathbb{H}^2 \times \mathbb{E}$. By their complicated nature, the last three geometric spaces are harder to understand and to work with than the Euclidean, hyperbolic, or spherical spaces. Thus, it is particularly important to provide interactive visualizations in order to get a feeling for these geometries. To this end, projects are currently underway (Coulon et al. 2020a; Coulon et al. 2020b; Novella, da Silva, and Velho 2020). However, the power of these visualizations is not limited to geometry education. A further engaging topic is the realm of science communication and optical illusions like the Penrose triangle, an impossible structure. That is to say, an impossible structure in three-dimensional Euclidean space, as it turns out that Penrose triangles do exist in $Nil$ geometry (Celińska-Kopczyńska and Kopczyński 2020).

Since the early work of Euclid, mathematicians have struggled with the existence of non-Euclidean geometries for two millennia. These spaces were unintuitive and too removed from the surrounding Euclidean world to lend themselves to distinct applications. The mathematics of Gauss, Bolyai, and Lobachevsky as well as Minkowski's and Einstein's physics finally provided a foundation that made use of curved spaces and their geometries. Compared to these long time-scales, Thurston's geometrization conjecture is still young. It remains to be seen what applications science will find for the geometries arising from it. Nevertheless, as for hyperbolic geometry, models and visualizations will play an important role. The ability to experience Thurston's geometries in an immersive way enables new research directions, new hypotheses, and new ways to view the material. Escher's woodcuts introduced generations of people to two-dimensional hyperbolic space. In much the same way, virtual reality is the perfect tool at the perfect time for the community to witness three-dimensional hyperbolic space.

---

[40] This fact goes by the name *uniformization theorem*. It has its own rich history, discussed in (Abikoff 1981).